\documentclass[12pt]{article}
\usepackage{amsmath,amssymb,amsbsy,amsfonts,amsthm,latexsym,amsopn,amstext,
amsxtra,euscript,amscd}

  \newtheorem{thm}{Theorem}%% [section]
  
  \newtheorem{lem}[thm]{Lemma}
  
  \theoremstyle{definition}
  
  \theoremstyle{remark}

\begin{document}

\def\squareforqed{\hbox{\rlap{$\sqcap$}$\sqcup$}}
\def\qed{\ifmmode\squareforqed\else{\unskip\nobreak\hfil
\penalty50\hskip1em\null\nobreak\hfil\squareforqed
\parfillskip=0pt\finalhyphendemerits=0\endgraf}\fi}

\def\cA{{\mathcal A}}
\def\cB{{\mathcal B}}
\def\cC{{\mathcal C}}
\def\cD{{\mathcal D}}
\def\cE{{\mathcal E}}
\def\cF{{\mathcal F}}
\def\cG{{\mathcal G}}
\def\cH{{\mathcal H}}
\def\cI{{\mathcal I}}
\def\cJ{{\mathcal J}}
\def\cK{{\mathcal K}}
\def\cL{{\mathcal L}}
\def\cM{{\mathcal M}}
\def\cN{{\mathcal N}}
\def\cO{{\mathcal O}}
\def\cP{{\mathcal P}}
\def\cQ{{\mathcal Q}}
\def\cR{{\mathcal R}}
\def\cS{{\mathcal S}}
\def\cT{{\mathcal T}}
\def\cU{{\mathcal U}}
\def\cV{{\mathcal V}}
\def\cW{{\mathcal W}}
\def\cX{{\mathcal X}}
\def\cY{{\mathcal Y}}
\def\cZ{{\mathcal Z}}

\def\fH{{\mathfrak H}}

\def \C {{\mathbb C}}
\def \F {{\mathbb F}}
\def \L {{\mathbb L}}
\def \K {{\mathbb K}}
\def \Q {{\mathbb Q}}
\def \Z {{\mathbb Z}}

\def\\{\cr}
\def\({\left(}
\def\){\right)}
\def\fl#1{\left\lfloor#1\right\rfloor}
\def\rf#1{\left\lceil#1\right\rceil}

\def \lcm{{\mathrm {lcm}}}
\def \llog{{\mathrm{llog~}}}
\def\e{\mathbf{e}}
\def\em{\e_m}
\def\el{\e_\ell}
\def\Res{\mathrm{Res}}

\def\ellax{\vec{a} \cdot \vec{x}}
\def\ellay{\vec{a} \cdot \vec{y}}

\def\mand{\qquad \text{and} \qquad}
\renewcommand{\vec}[1]{\mathbf{#1}}

\newcommand{\comm}[1]{\marginpar{%
\vskip-\baselineskip %raise the marginpar a bit
\raggedright\footnotesize
\itshape\hrule\smallskip#1\par\smallskip\hrule}}

\title{On the Greatest Common Divisor
% and
%the Least Common Multiple
of Shifted Sets}

\author{
%{\sc Etienne Fouvry}  \\
%{D$\acute{\textrm{e}}$partment de Math$\acute{\textrm{e}}$matiques,
%Universit$\acute{\textrm{e}}$ Paris-Sud}\\
%{91405 Orsay Cedex, France}\\
%{\tt etienne.fouvry@math.u-psud.fr}
%\and
%{\sc Joachim von zur Gathen}  \\
%{B-IT, Universit$\ddot{\textrm{a}}$t Bonn}\\
%{53113 Bonn, Germany}\\
%{\tt gathen@bit.uni-bonn.de}
%\and
{\sc Randell Heyman and Igor E. Shparlinski}\\
{School of Mathematics and Statistics}\\
{University of New South Wales} \\
{Sydney, NSW 2052, Australia}\\
{\tt \{randell,igor.shparlinski\}@unsw.edu.au}
%\and
%{\sc Igor E. Shparlinski}  \\
%{School of Mathematics and Statistics, University of New South Wales} \\
%{Sydney, NSW 2033, Australia}\\
%{\tt igor.shparlinski@unsw.edu.au}
}

\date{\today}

\maketitle

\begin{abstract}
Given a set  of $n$ positive  integers $\{a_1, \ldots, a_n\}$
and an integer parameter $H$
we study small additive shift of its elements by
integers $h_i$ with $|h_i| \le H$, $i =1, \ldots, n$, such that the
greatest common divisor of
$a_1+h_1, \ldots, a_n+h_n$ is very different from that of $a_1, \ldots, a_n$.
We also consider a similar problem for the least common multiple.
\end{abstract}

%\paragraph{Mathematical Subject Classification (2000):} Primary ???; Secondary
%???

%\paragraph{Keywords:} GCD, LCM,  circulant graphs, approximate GCD

\section{Introduction}

Let $\vec{a} = (a_1, \ldots, a_n) \in \Z^n$ be a nonzero vector.
The {\it approximate common divisor problem\/}, introduced by Howgrave-Graham~\cite{How}
for $n=2$, can generally be described as follows. Suppose we are given  two bounds $D> h\ge 1$.
Assuming that for some $h_i$ with $|h_i| \le H$, $i =1, \ldots, n$,
we have
\begin{equation}
\label{eq:gcd D}
\gcd(a_1+h_1, \ldots, a_n+h_n) > D,
\end{equation}
the task is to determine the shifts $h_1, \ldots, h_n$.
If it is also requested that $h_1=0$ then we refer to the problem as the {\it partial
approximate common divisor problem\/} (certainly in this case the task is to
find the shifts faster than via complete factorisation of $a_1\ne 0$).

This problem has a strong cryptographic motivation as it is related to
some attacks on the RSA and some other cryptosystems, see~\cite{ChNg,CoHe,How,SaMa} and references therein
for various algorithms and applications.
In particular,  much of the current motivation for studying approximate common divisor problems
stems from the search for efficient
and reliable {\it fully homomorphic encryption\/}, that is, encryption that allows arithmetic  operations
on encrypted data, see~\cite{vDGHV,Gen,Mic}.

Here we consider a dual question and show that for any
$\vec{a} = (a_1, \ldots, a_n) \in \Z^n$, there are shifts $|h_i| \le H$, $i =1, \ldots, n$,
for which~\eqref{eq:gcd D} holds with a relatively large value of $D$. Throughout we  use
$\gcd(\vec{x})$  to mean $\gcd(x_1,\ldots,x_n)$
 for any $\vec{x} \in \Z^n$.

We also denote the height of $\vec{x}$ with
$\fH(\vec{x})=\max\{|x_1|, \ldots, |x_n|\}$.

The implied constants in the symbols `$O$',  `$\ll$'
and `$\gg$'  may occasionally,
where obvious, depend  on the integer parameter $n$
and the real positive parameter $\varepsilon$,
and are absolute otherwise.
We recall that the notations $U = O(V)$,  $U \ll V$  and $V \gg U$ are all
equivalent to the assertion that the inequality $|U|\le c|V|$ holds for some
constant $c>0$.

Our treatment of this question is based on some results of Baker and Harman~\cite{BaHa2}
(see also~\cite{BaHa1}). For an integer $n \ge 1$ and real
positive $\varepsilon<1$, we define
$\kappa(n, \varepsilon)$ as the solution $\kappa>0$ to the equation
\begin{equation}
\label{eq:kappa}
\frac{n(\varepsilon \kappa - 1)}{n-1} = \frac{1}{2^{2 +\max\{1,\kappa\}} - 4}.
\end{equation}
The solution is unique as the
left hand side of~\eqref{eq:kappa} is monotonically increasing (as a
function of $\kappa$) from $-n/(n-1)$ to $+\infty$ on
$[0, \infty)$ while the right hand
side of~\eqref{eq:kappa} is positive and monotonically non-increasing.

We also set
$$
\vartheta(n, \varepsilon)
= \frac{1}{(n-1)} \(1 - \frac{1}{\varepsilon \kappa(n, \varepsilon)}\).
$$
It  easy to see from~\eqref{eq:kappa} that $\varepsilon \kappa(n, \varepsilon) < 1$,
so $\vartheta(n, \varepsilon) >0$.

\begin{thm}\label{thm:gcd large}
\label{ques:gcd shift}  For any
vector $\vec{a} \in \Z^n$, any real positive $\varepsilon <1$
and
$$
H\ge \fH(\vec{a})^\varepsilon
$$
there exists a vector $\vec{h} = (h_1, \ldots, h_n) \in \Z^n$
of height
$$\fH(\vec{h})  \le  H
$$
such that
$$
\gcd(\vec{a+h}) \gg \fH(\vec{h})  H^{\vartheta(n, \varepsilon)}.
$$
\end{thm}

Next we are  interested in asking for which $\vec{h}$ the shifted set is pairwise
coprime.

For $\vec{a} \in \Z^n$ we denote by $L(\vec{a})$ the smallest $H$ such that there is a
$\vec{h} \in \Z^n$ with $\fH(\vec{h})=H$ such that
$$
\gcd(a_i+h_i, a_j + h_j) =1, \qquad 1 \le i < j \le n.
$$
For $n=2$, and thus $\vec{a}=(a_1,a_2) \in \Z^2$, 
Erd{\H o}s~\cite[Equation~(3)]{Erd} has given the bound
$$
L(\vec{a})\ll
\frac{\log \min\{|a_1|, |a_2|\}}{\log\log \min\{|a_1|, |a_2|\}}.
$$
However the method of~\cite{Erd} does not seem to generalise 
to $n \ge 3$.

\begin{thm}\label{thm:gcd 1}
For  an arbitrary $\vec{a} \in \Z^n$ we have
$$L(\vec{a})\ll\log^2 \fH(\vec{a}).$$
\end{thm}

Note in fact our argument allows to replace $\fH(\vec{a})$ with a smaller
qunatity
$$
\fH^*(\vec{a}) = \min_{1 \le i \le n} 
\max_{\substack{1 \le j \le n\\ i\ne j}} |a_i|.
$$
For $\vec{a} \in \Z^n$ we denote by $\ell(\vec{a})$ the smallest $H$ such that there is a vector
$\vec{h} \in \Z^n$ with $\fH(\vec{h})=H$ and
$$
\gcd(a_1+h_1, \ldots, a_n + h_n) =1.
$$

A very simple argument, based on the Chinese Remainder Theorem,
implies the following result, which generalises~\cite[Equation~(2)]{Erd}.

\begin{thm}\label{thm:crt}
For infinitely many $\vec{a} \in \Z^n$ we have
$$\ell(\vec{a})\gg
\(\frac{\log \fH(\vec{a})}{\log\log  \fH(\vec{a})}\)^{1/n} .
$$
\end{thm}

Note that Theorem~\ref{thm:crt} is essentially an explicit version of
a result of Huck and   Pleasants~\cite{HuPl}.

It is clear that for non-zero vector  $\vec{a} \in \Z^n$
and arbitrary vectors $\vec{x} , \vec{y}  \in \Z^n$
we have
$$\gcd(a_1, \ldots, a_m)\mid \gcd(\ellax,\ellay),
$$
where
$$
\ellax = \sum_{i=1}^n a_i x_i \mand \ellay = \sum_{i=1}^n a_i y_i.
$$
Let   $R(\vec{a}, h)$ be the
number of vectors $\vec{x}, \vec{y}\in \Z^n$ with
positive components and  of height
$\fH(\vec{x}), \fH(\vec{y})\le h$ for which
\begin{equation}
\label{GCD=gcd}
\gcd(a_1, \ldots, a_m)= \gcd(\ellax,\ellay).
\end{equation}

By~\cite[Theorem~3]{vzGSh} we have
$$
|R(\vec{a}, h) - \zeta(2)^{-1}h^{2n} | \le
h^{2n-1/n}(h \fH(\vec{a}) )^{o(1)},
$$
where $\zeta(s)$ is the Riemann zeta function.

\begin{thm}
\label{thm:gcd lin form} Let $n \ge 2$ and let
$\vec{a} \in \Z^n$.
Then, for  $\max\{h,\fH(\vec{a})\}\to \infty$,
$$
|R(\vec{a}, h) - \zeta(2)^{-1}h^{2n} | \le
h^{2n-n/(n^2 - n +1)}(h \fH(\vec{a}))^{o(1)}.
$$
\end{thm}

\section{Proof of Theorem~\ref{thm:gcd large}}

We use the following~\cite[Theorem~1]{BaHa2}, see also~\cite[Equation~(2.1)]{BaHa2}
that gives an explicit formula for constant $\gamma(K)$ below.

\begin{lem}\label{lem:BH}  Suppose that for some fixed $K>0$
and some sufficiently large real positive $Q$ and $R$
we have
$$
\(\sum_{i=1}^n a_i^2\)^{1/2} \le  R^K
$$
and
$$
C_1(K,n) \le Q \le R^{\gamma(K)},
$$
where
$$
\gamma(K) = \frac{1}{2^{2+\max\{1,K\}} - 4}.
$$
Let $\psi_1,\ldots,\psi_n$ be positive integers with
$$\psi_i \le c_2(K,n)  (\log Q)^{-n}, \quad i=1,\ldots,n,$$

and
$$\psi_1\cdots\psi_n=Q^{-1}.$$
Then $$\left\|\frac{a_i}{r}\right\|\le \psi_i,\quad i=1,\ldots,n,$$
$$R \le r \le 2QR.$$
where $C_1(K,n)$ and $c_2(K,n)$ depend at most on $K$ and $n$.
\end{lem}

To prove Theorem~\ref{thm:gcd large}, we choose some parameters $Q$ and $R$ that
satisfy Lemma~\ref{lem:BH} with $K = \kappa(n,\varepsilon)$,
where $ \kappa(n,\varepsilon)$ is given  by~\eqref{eq:kappa},
and then we set $\psi_i=Q^{-1/n}$, $i=1, \ldots, n$.
Then by Lemma~\ref{lem:BH}, there exist an integer $r$
with $R \le r \le 2QR$ such that
$$
\left\|\frac{a_i}{r}\right\| \le Q^{-1/n}, \qquad i = 1, \ldots, n,
$$
where $\|\xi\|$ is distance between a real $\xi$ and
the closest integer.
So for some integers $h_i$ with $|h_i| \le r Q^{-1/n}$ we have
$$
a_i + h_i \equiv 0 \pmod r, \qquad i = 1, \ldots, n.
$$

Suppose that for some constant $A > 0$
we choose $R$ such that for $Q = (0.5)^{n/(n-1)}  A^{-1} R^{\gamma(K)}$, we have
\begin{equation}
\label{RfromH}
2Q^{1-1/n}R =H.
\end{equation}
Then
$$
R = A^{(n-1)/(n\gamma(K)+n)} H^{n/(n\gamma(K)-\gamma(K)+n)}.
$$
Then,
taking $A$ to satisfy
$$
A^{(n-1)/(n\gamma(K)+n)} = n^{1/2K}
 $$
 due to our choice of $K = \kappa(n,\varepsilon)$, we have
\begin{equation}
\label{eq:R H}
R  =n^{1/2K} H^{n/(n\gamma(K)-\gamma(K)+n)} = n^{1/2K}H^{1/\varepsilon K}.
\end{equation}
Hence for $\vec{h} = (h_1, \ldots, h_n)$ we have
$$
\fH(\vec{h}) \le r Q^{-1/n} \le 2Q^{1-1/n}R =H
$$
and
\begin{equation}
\label{eq:gcd Q}
\gcd(\vec{a+h}) \ge r\ge  \fH(\vec{h}) Q^{1/n}. 
\end{equation}

Using~\eqref{eq:R H}, we derive
$$
\(\sum_{i=1}^n a_i^2\)^{1/2}  \le n^{1/2} H^{1/\varepsilon} = R^{K}.
$$
Thus Lemma~\ref{lem:BH} indeed applies.
We also have
\begin{equation}
\label{eq:Q lower}
Q^{1/n} \gg  R^{\gamma(K)/n} \gg H^{\gamma(K) /\varepsilon n K}.
\end{equation}
We now see from~\eqref{eq:kappa} that
$$
\frac{\gamma(K)}{\varepsilon n K} =
\frac{\varepsilon K-1}{\varepsilon (n-1) K},
$$
which together with~\eqref{eq:gcd Q} and~\eqref{eq:Q lower} completes the proof.

\section{Proof of Theorem~\ref{thm:gcd 1}}

We recall  the following well-known result of Iwaniec~\cite{Iw}
on the {\it Jacobsthal problem\/}.
For a  given $r$, let $C(r)$ be the maximal length of a sequence of consecutive integers, each divisible by one of $r$ arbitrarily chosen primes. Then Iwaniec~\cite{Iw} gives the following
bound on $C(r)$:

\begin{lem}\label{Iwaneic} For a given $r>1$ we have,
$$C(r) \ll (r \log r)^2.$$
\end{lem}

We are now ready to prove Theorem~\ref{thm:gcd 1}.

We now set $h_1=0$ and chose $h_i$, $i =2, \ldots, n$
as the smallest non-negative integer with
$$\gcd\(\prod_{j=1}^{i-1} (a_{j}+h_{j}),a_i+h_i\)=1.
$$

We show that if $n$ is a positive integer and $a=\fH(\vec{a})$
then
\begin{equation}\label{statement}
\fH(\vec{h})\ll\log^2 a.
\end{equation}

For $n=2$ we note that $a_1$ has $\omega(a_1)$ distinct prime factors,
where $\omega(a)$ is the number of distinct
prime divisors of an integer  $a\ge 1$.

So, by Lemma~\ref{Iwaneic},
$$C(\omega(a_1)) \ll (\omega(a_1)\log(\omega(a_1))^2\ll \log^2 a_1 = \log^2 a$$
for all $a_1$, and from the trivial bound $\omega(k)! \le k$ and
the Stirling formula we have
 $$\omega(k)  \ll \frac{\log k}{\log(2+ \log k)}$$
 for any integer $k \ge 1$.
Now a straight forward inductive argument, after simple calculations,
implies~\eqref{statement} and concludes the proof.

\section{Proof of Theorem \ref{thm:crt}}

Let us choose a sufficiently large parameter $H$ and the first $(2H+1)^n$
 primes  $p_{i_1,\ldots,i_n} > H$ for $-H \le i_1,\ldots,i_n \le H$.

For each $k=1, \ldots, n$ we define $a_k$ as the smallest
positive integer
with
$$
a_k  \equiv i_k \pmod {p_{i_1,\ldots,i_n}}, \qquad -H \le i_1,\ldots,i_n \le H.
$$
Set $\vec{a} = (a_1, \ldots,a_n)$. Clearly, for any $\vec{h} \in \Z^n$
with $\fH(\vec{h})\le H$, we have
$$
p_{h_1,\ldots,h_n} \mid \gcd(a_1+h_1, \ldots, a_n + h_n) .
$$
This implies that $\ell(\vec{a}) \ge H$.

It remains to estimate $\fH(\vec{a})$.
Clearly, we have $p_{i_1,\ldots,i_n} \ll H^n \log H$ for $-H \le i_1,\ldots,i_n \le H$.
Therefore,
$$
\fH(\vec{a}) \le \prod_{-H \le i_1,\ldots,i_n \le H} p_{i_1,\ldots,i_n}
= \exp(O(H^n \log H)) = \exp(O(\ell(\vec{a})^n \log \ell(\vec{a}))),
$$
which completes the proof.

\section{Proof of Theorem~\ref{thm:gcd lin form}}

Clearly, it is enough to
consider the case where $\gcd(a_1, \ldots, a_n) = 1$.

We can certainly assume that $n \le \log h$ for otherwise the bound
is trivial.

Let  $\mu$ denote the  M\"obius function, that is $\mu(1)= 1$, $\mu(d) = 0$ if $d \ge 2$ is not squarefree,
and $\mu(d) = (-1)^{\omega(d)}$  otherwise, where $\omega(d)$, as before,
is the number of
prime divisors of an integer  $d\ge 1$.

As in the proof of~\cite[Theorem~3]{vzGSh}, by the inclusion exclusion principle we have
$$R(\vec{a}, h) = \sum_{d\geq 1} \mu(d) U_d(\vec{a}, h) ^2,
$$
where for an integer  $d \ge 1$, we denote by  $U_d(\vec{a}, h)$
the number of vectors $\vec{x} \in \Z^n$ with
positive components and  of height
$\fH(\vec{x})\le h$ for which  $d \mid \ellax$.

We now recall from~\cite{vzGSh} some estimates on $U_d(\vec{a}, h)$.

More precisely, for $1 \le d \le 2h/3n$ we have
\begin{equation}
\label{eq:Asymp U}
\left|U_d(\vec{a}, h)^{2}-\frac{h^{2n}}{d^{2}}\right|\le 8n d^{-1} h^{2n-1}.
\end{equation}
see~\cite[Equation~(8)]{vzGSh}.
The proof of~\eqref{eq:Asymp U} also relies on the bound
\begin{equation}
\label{eq:Bound U1}
U_d(\vec{a}, h) \le d^{n-1}\(h/d   + 1\)^n .
\end{equation}
that holds for any integer $d \ge 1$.

Furthermore,
for any squarefree $d \ge 1$ we also have  the bound
\begin{equation}
\label{eq:Bound U2}
U_d(\vec{a}, h) \le h^{n-1}\(h  d^{-1/n} + 1\) .
\end{equation}
see~\cite[Equation~(10)]{vzGSh}.

Therefore, choosing some parameter $D$,  we write
\begin{equation}
\label{eq:R expand}
R(\vec{a}, h) =  M + O(\Delta_1 + \Delta_2)
\end{equation}
where
\begin{align*}
M &=   \sum_{d \le 2h/3n } \mu(d)  U_d(\vec{a}, h)^2,\\
\Delta_1 &= \sum_{2h/3n < d \le D} \mu(d)  U_d(\vec{a}, h)^2,\\
\Delta_2 &= \sum_{d >D} \mu(d)  U_d(\vec{a}, h)^2.
\end{align*}

Using~\eqref{eq:Asymp U}, we derive
\begin{align*}
M &=
\sum_{d \le 2h/3n} \mu(d) \(\frac{h^{2n}}{d^2}
+ O\(h^{2n-1} d^{-1}\) \) \\
& = h^{2n}\sum_{d \le 2h/3n} \frac{\mu(d)}{d^2}  +  O\(h^{2n-1}\log h\).
\end{align*}
Since
$$
\sum_{d \le 2h/3n} \frac{\mu(d) }{d^2}
= \sum_{d= 1}^\infty \frac {\mu(d)
}{d^2} +  O\(D^{-1}\)  = \zeta(2)^{-1}  +    O\(D^{-1}\),
$$
see~\cite[Theorem~287]{HaWr},
we derive
\begin{equation}
\label{eq:M asymp}
M = h^{2n} \zeta(2)^{-1}  +    O\(h^{2n-1}\log h\).
\end{equation}

To  estimate $\Delta_1$ we apply the bound~\eqref{eq:Bound U1},
which for $d \ge 2h/3n$ can be simplified as $U_d(\vec{a}, h)  = O(d^{n-1})$.
Therefore,
\begin{equation}
\label{eq:Delta1 prelim}
\Delta_1 \ll  \sum_{2h/3n < d \le D}  d^{n-1} U_d(\vec{a}, h)
\le D^{2n-1} \sum_{2h/3n < d \le D} U_d(\vec{a}, h) .
\end{equation}
Using the same argument as the proof of~\cite[Theorem~3]{vzGSh},
based on a bound of the divisor function $\tau(k)$, we
obtain
\begin{equation}
\begin{split}
\label{eq:Sum U}
        \sum_{d > D} U_d(\vec{a}, h)& =
\sum_{ d > D} \sum_{\substack{\fH(\vec{x}) \le h\\
d \mid \ellax}} 1\\
&  =  \sum_{h(\vec{x}) \le h}   \sum_{\substack{d >  D \\ d \mid \ellax }} 1
\le
\sum_{h(\vec{x}) \le h}  \tau(\ellax ) \le h^n (h \fH(\vec{a}))^{o(1)},
\end{split}
\end{equation}
where $\vec{x}$ runs through integral vectors with positive
components.  Hence, we see that~\eqref{eq:Delta1 prelim} yields
the estimate
\begin{equation}
\label{eq:Delta1 bound}
\Delta_1 \ll D^{n-1}h^n (h \fH(\vec{a}))^{o(1)}.
\end{equation}

Finally, to estimate $\Delta_2$ we apply the bound~\eqref{eq:Bound U2}
and, as before derive
\begin{equation}
\label{eq:Delta2 bound}
\Delta_2 \ll h^{n-1}\(h D^{-1/n} + 1\) \sum_{ d > D}
U_d(\vec{a}, h)  \le h^{2n-1}\(h D^{-1/n} + 1\) (h  \fH(\vec{a}))^{o(1)}.
\end{equation}

Substituting the bounds~\eqref{eq:M asymp}, \eqref{eq:Delta1 bound}
and~\eqref{eq:Delta2 bound} into~\eqref{eq:R expand}, we obtain
$$
R(\vec{a}, h) =    h^{2n} \zeta(2)^{-1}  +
O\(\(h^{2n-1}+D^{n-1}h^n + h^{2n}D^{-1/n}\) (h \fH(\vec{a}))^{o(1)} \).
$$
Now, choosing
$$
D = h^{n^2/(n^2 - n +1)},
$$
we conclude the proof.

\section{Comments}

We remark that it is also interesting to study
analogous questions for polynomials with integer
coefficients or over finite fields, see~\cite{EGB,Na,vzGMS}
for some polynomial versions of the
approximate common divisor problem. Some of out techniques can
be extended to this case, however some important ingredients,
such as the results of Baker and Harman~\cite{BaHa1,BaHa2}
are missing.

\section*{Acknowledgment}

The authors are grateful to Etienne Fouvry and Michel Laurent for 
very useful discussions. 

This work was supported in part by the ARC Grant DP130100237.

\end{document}